\documentclass{article}

\usepackage{amssymb}
\usepackage{amsmath}
\usepackage{amsthm}

\newtheorem{theorem}{Theorem} \newtheorem{lemma}{Lemma}[section]
\newtheorem{propo}{Proposition}[section]

  \newcommand{\ep}{\varepsilon}
  
 \newcommand{\e}{\ep} 
\newcommand{\F}{\mathbb{F}} 
\newcommand{\G} {\Gamma}

\newcommand{\Ker}{\mbox{Ker}}

\newcommand{\cA}{\mathcal{A}} \newcommand{\g}{\gamma}
\newcommand{\Sub}{\mbox{Sub}}
\newcommand{\Stab}{\mbox{Stab}}
\newcommand{\vi}{\vskip 0.1in \noindent}
\newcommand{\roo}{\mbox{root}}
\newcommand{\Forg}{\mbox{F}}
\newcommand{\dist}{\mbox{dist}}
\newcommand{\diam}{\mbox{diam}}\newcommand{\Cay}{\mbox{Cay}}
\newcommand{\sk}{S^{K,Q}_\G}
\title{On uniformly recurrent subgroups of finitely generated groups\footnote{AMS
Subject Classification: 37B05, 20E99} }

\author{G\'abor Elek}

\begin{document}
\maketitle

\begin{abstract}
We prove that if $G$ is a finitely generated group and $Z$ is a uniformly
recurrent subgroup of $G$ then there exists a minimal system $(X,G)$ with
$Z$ as its stability system. This answers a query of Glasner and Weiss
\cite{GW} in the case of finitely generated groups. Using the same method
(introduced by Alon,  Grytczuk, Haluszczak and Riordan \cite{AGHR}) we will
prove that finitely generated sofic groups have free Bernoulli-subshifts
admitting an invariant probability measure.
\end{abstract}

\textbf{Keywords.} uniformly recurrent subgroups, sofic groups
\section{Introduction}
Let $\G$ be a countable group and $\Sub(\G)$ be the compact space of
all subgroups of $\G$. The group $\G$ acts on $\Sub(\G)$ by conjugation.
{\it Uniformly recurrent subgroups}\,(URS) were defined by 
Glasner and Weiss \cite{GW}
as closed, invariant subsets $Z\subset \Sub(\G)$ such that the action
of $\G$ on $Z$ is minimal (every orbits are dense).
Now let $(X,\G,\alpha)$ be a $\G$-system (that is, $X$ is a compact metric
space and $\G\to \mbox{Homeo}(X)$ is a homomorphism).
For each point $x\in X$ one can define the topological stabilizer subgroup 
$\Stab_\alpha^0(x)$ by
$$\Stab_\alpha^0(x)=\{\g\in \G\,\mid\, \g \,\mbox {fixes some neighborhood of $x$}\}\,.$$

\noindent
Let us consider the  $\G$-invariant subset $X^0\subseteq X$ such that
$x\in X^0$ if and only if $\Stab_\alpha(x)=\Stab_\alpha^0(x)\,.$ Then
$X^0$ is a dense $G_\delta$-set and  we have a $\G$-equivariant map 
$S_\alpha:X^0\to \Sub(\G)$ such that
if $y\in X^0$ then $S_{\alpha}(y)=\Stab_\alpha(y)$. The closure of the
invariant subset 
$S_{\alpha} (X^0)\subset \Sub(\G)$ is
called the {\it stability system} of $(X,\G,\alpha)$ 
(see also \cite{Kaw},\cite{BB}).
If the action is minimal, then the stability system of $(X,\G,\alpha)$
is an URS.
Glasner and Weiss proved (Proposition 6.1,\cite{GW}) that
for every URS $Z\subset \Sub(G)$ there exists a topologically transitive
(that is there is a dense orbit) system $(X,\G,\alpha)$ with $Z$ as its
stability system.
They asked (Problem 6.2., \cite{GW}), whether for any URS $Z$
there exists a 
minimal system $(X,\G,\alpha)$ with $Z$ as its stability system. Recently,
Kawabe \cite{Kaw} gave an affirmative answer for this question in the case
of amenable groups. We will prove the following result.
\begin{theorem}\label{tetel1}
If $\G$ is a finitely generated group and $Z\subset \Sub(\G)$ is an URS, then
there exists a minimal system $(X,\G,\alpha)$ with $Z$ as its stability
system.
\end{theorem}
\noindent
In the proof we will use the Lov\'asz Local Lemma technique of 
Alon, Grytczuk, Haluszczak and Riordan
 \cite{AGHR} to construct a minimal action
on the space of rooted colored $\G$-Schreier graphs. 
This approach has already
been used to construct free $\G$-Bernoulli subshifts by Aubrun, Barbieri and Thomass\'e 
\cite{ABT} . The other result of the paper is about
free $\Gamma$-Bernoulli-subshifts, that is
closed $\Gamma$-invariant subsets $M$ of $K^\Gamma$, where
$K$ is some finite alphabet and the action
of $\Gamma$ on $M$ is free.  For a long time all finitely generated groups that
had been known to have free Bernoulli-subshifts were residually-finite.
Then Dranishnikov and Schroeder \cite{DS} constructed a free Bernoulli-subshift
for any torsion-free hyperbolic group. Somewhat later Gao, Jackson and Seward proved that any countable group has free Bernoulli-subshifts \cite{GJS1}, 
\cite{GJS2}. 
On the other hand, Hjorth and Molberg \cite{HM} proved that for any countable group $\Gamma$ there exists a free continuous action of $\Gamma$
on a Cantor set admitting an invariant measure. It seems that
so far all groups $\Gamma$ for which free Bernoulli-shifts with
an invariant probability measure proved to exist were either residually-finite (Toeplitz-shubshifts) or amenable (when the existence of invariant measure is
obvious). Using the  coloring technique of Alon, Grytczuk, Haluszczak and Riordan we prove
a combination of these results for finitely generated sofic groups.
\begin{theorem}\label{tetel2}
Let $\G$ be a finitely generated sofic group. Then there exists a free Bernoulli-subshift for $\Gamma$.
\end{theorem}

\section{The space of colored rooted $\G$-Schreier graphs}
Let $\Gamma$ be a finitely generated group with a minimal symmetric
generating system $Q=\{\gamma_i\}^n_{i=1}$. 
Let $H\in\Sub(\Gamma).$ Then the Schreier graph of $H$, $S^Q_\G(H)$ is 
constructed as follows.
\begin{itemize}
\item The vertex set of $S^Q_\G(H)$, $V(S^Q_\G(H))=\Gamma/H$ (that is
$\Gamma$ acts on the vertex set of $S^Q_\G(H)$ on the left).
\item The vertices corresponding to the cosets $aH$ and $bH$ are connected
by a directed edge labeled by the generator $\gamma_i$ if $\gamma_iaH=bH$.
\end{itemize}
\vi
The coset class of $H$ is called the {\it root} of the graph  $S^Q_\G(H)$.
We will consider the usual shortest path distance on  $S^Q_\G(H)$ and
denote the ball of radius $r$ around the root $H$ by $B_r(S^Q_\G(H), H)\,.$
Note that  $B_r(S^Q_\G(H), H)$ is a rooted edge-labeled graph.
The space of all Schreier graphs $S^Q_\Gamma$ is a compact metric
space, where
$$d_{S^Q_\G}(S^Q_\G(H_1), S^Q_\G(H_2))=2^{-r}\,,$$
if $r$ is the largest integer for which the $r$-balls
$B_r(S^Q_\G(H_1), H_1)$ and \\ $B_r(S^Q_\G(H_2), H_2)$ are rooted-labeled isomorphic.
Clearly, $s:\Sub(\Gamma)\to S^Q_\Gamma$, $s(H)=S^Q_\G(H)$ is a 
homeomorphism commuting
with the $\G$-actions. Note that if $\gamma\in\Gamma$ and $H\in\Sub(\Gamma)$,
then 
$$\gamma(S^Q_\Gamma(H))=S^Q_\Gamma(\gamma H \gamma^{-1})\,,$$
where the underlying labeled graphs of $S^Q_\Gamma(H)$ and 
$S^Q_\Gamma(\gamma H \gamma^{-1})$ are isomorphic. The graph 
$S^Q_\Gamma(\gamma H \gamma^{-1})$ can be regarded as the same graph
as $S^Q_\Gamma(H)$ with the new root $\gamma(\roo(S^Q_\G(H)))\,.$
We will use the root-change picture of the $\Gamma$-action on
$S^Q_\G$ later in the paper.
\vi
Now let $K$ be a finite alphabet. A rooted $K$-colored
Schreier graph is a rooted Schreier graph $S^Q_H$ equipped
with a vertex-coloring $c:\Gamma/H\to K$. 
Let $\sk$ be the set of all rooted $K$-colored Schreier-graphs. Again,
we have a compact, metric topology on $\sk$:
$$d_{\sk}(S,T)=2^{-r}\,,$$
\vi
if $r$ is the largest integer such that the $r$-balls around the
roots of the graphs $S$ and $T$ are rooted-colored-labeled isomorphic.
We define $d_{\sk}(S,T)=2$ if the $1$-balls around the roots are nonisomorphic
and even the colors of the roots are different.
Again, $\Gamma$ acts on the compact space $\sk$ by the root-changing map.
Hence, we have a natural color-forgetting map $\Forg:\sk\to S^Q_\Gamma$
that commutes with the $\Gamma$-actions. Notice that if
a sequence $\{S_n\}^\infty_{n=1}\subset \sk$ converges to $S\in\sk$, then 
for any $r\geq 1$ there exists some integer $N_r\geq 1$ such that
if $n\geq N_r$ then the $r$-balls around the roots of the graph $S_n$ and
the graph $S$ are rooted-colored-labeled isomorphic. 
Let $H\in\Sub(\Gamma)$ and $c:\G/H\to K$ be a vertex coloring
that defines the element  $S_{H,c}\in \sk$. Then of course,
$\gamma(S_{H,c})=S_{H,c}$ if $\gamma\in H$. On the other hand, if
$\gamma(S_{H,c})=S_{H,c}$ and $\gamma\notin H$ then we have the
following lemma that is immediately follows from the definitions of the $\G$-actions.
\begin{lemma}\label{trivi}
Let $\gamma\notin H$ and $\gamma(S_{H,c})=S_{H,c}$. Then there exists
a colored-labeled graph-automorphism of the $K$-colored
labeled graph $S_{H,c}$ moving
the vertex representing  $H$ to the vertex representing $\gamma(H)\neq H$.
\end{lemma}
\vi
Note that we have a continuous $\G$-equivariant map $\pi:\sk\to\Sub(\Gamma)$,
where $\pi(t)=s^{-1}\circ \Forg(t).$
Let $Z$ be an URS of $\Gamma$. Let $H\in Z$ and let $t\in \sk$ be
corresponding to a vertex coloring of the
Schreier graph $S^Q_\Gamma(H)$. We say that the element
$t\in \sk$ is $Z$-proper if $\Stab_\alpha(t)=H$, where $\alpha$ is the
right action of $\Gamma$ on $\sk$. Note that if $H\in Z$ and $t$ is representing
the Schreier graph $S^Q_H$, then by Lemma \ref{trivi}, $t$ is $Z$-proper if and only
if there is no non-trivial colored-labeled automorphism of $t$.
\begin{propo}\label{propo1}
Let $Y\subset \sk$ be a closed $\G$-invariant subset consisting
of $Z$-proper elements. Let
$(M,\G,\alpha)\subset (Y,G,\alpha)$ be a minimal $\Gamma$-subsystem.
Then for any $m\in M$, $\Stab_\alpha^0(m)=\Stab_\alpha(m)\in Z$. Also,
$\pi(M)=Z$.
\end{propo}
\proof
Let $h\in \Stab_\alpha(m)$. Then $h\in Z$, that is, $h$ fixes the root of $m$.
Therefore, $h$ fixes the root of $m'$ provided that $d_{\sk}(m,m')$ is small
enough.  Thus, $h\in \Stab_\alpha^0(m)$. Since $\pi$ is a
$\Gamma$-equivariant continuous map and $M$ is a closed $\Gamma$-invariant
subset, $\pi(M)=Z$. \qed
\section{The proof of Theorem \ref{tetel1}}
Let $Z$ be an URS of $\G$.
By Proposition \ref{propo1}, it is enough to construct a closed
$\Gamma$-invariant subset $Y\subset \sk$ for some alphabet $K$ such
that all the elements of $Y$ are $Z$-proper.
Let $H\in Z$ and consider the Schreier graph $S=S^Q_\G(H)$. Following \cite{ABT} and
\cite{AGHR} we call a coloring $c:\Gamma/H\to K$ nonrepetitive if
for any path $(x_1,x_2,\dots, x_{2n})$ in $S$ there exists some
$1\leq i\leq n$ such that
$c(x_i)\neq c(x_{n+i})\,.$ We call all the other colorings repetitive.
\begin{theorem}\label{alon}[Theorem 1 \cite{AGHR}]
For any $d\geq 1$ there exists a constant $C(d)>0$ such that any 
graph $G$ (finite or infinite) with vertex degree bound $d$ has
a nonrepetitive coloring with an alphabet $K$, provided that  $|K|\geq C(d)$.
\end{theorem}
\proof Since the proof in \cite{AGHR} is about edge-colorings and
the proof in \cite{ABT} is in slightly different setting, for completeness
we give a proof using  Lov\'asz's  Local Lemma, that closely follows the
proof in \cite{AGHR}. Now, let us state the Local Lemma.
\begin{theorem}[The Local Lemma]
Let $X$ be a finite set and $\Pr$ be a probability distribution on
the subsets of $X$. For $1\leq i\leq r$ let $\cA_i$ be a set of 
events, where an ``event'' is just a subset of $X$. Suppose that for
all $A\in \cA_i$, $\Pr(A_i)=p_i$.
 Let $\cA=\cup^r_{i=1} \cA_i$. Suppose that there are real numbers
$0\leq a_1,a_2,\dots,a_r<1$  and $\Delta_{ij}\geq 0$, $i,j=1,2,\dots, r$
such that the following conditions hold:
\begin{itemize}
\item for any event $A\in  \cA_i$ there exists a set $D_A\subset \cA$ with
$|D_A\cap \cA_j|\leq \Delta_{ij}$ for all $1\leq j \leq r$ such that
$A$ is independent of $\cA\backslash (D_A\cup \{A\}),$
\item $p_i\leq a_i\prod^r_{j=1}(1-a_j)^{\Delta_{ij}}$ for all $1\leq i \leq r\,.$
\end{itemize}
Then $\Pr(\cap_{A\in\cA} \overline{A})>0$. 
\end{theorem}
\vi
Let $G$ be a finite graph with maximum degree $d$. It is enough to
prove our theorem for finite graphs. Indeed, if $G'$ is a connected
infinite graph with vertex degree bound $d$, then for each ball around
a given vertex $p$ we have a nonrepetitive coloring. Picking a pointwise
convergent subsequence of the colorings we obtain a nonrepetitive 
coloring of our
infinite graph $G'$. 

\noindent
Let $C$ be a large enough number, its exact value will
be given later. Let $X$ be the set of all random $\{1,2,\dots,C\}$-colorings
of $G$. Let $r=\diam(G)$ and for $1\leq i \leq r$ and for any path $P$ of 
length $2i-1$ 
let $A(P)$ be the event that $P$ is repetitive. Set
$$\cA_i=\{A(P):\, \mbox{$P$ is a path of length $2i-1$ in $G$}\}\,.$$
\vi
Then $p_i=C^{-i}$. The number of paths of length $2j-1$ that intersects
a given path of length $2i-1$ is less or equal than $4ij d^{2j}$.
So, we can set $\Delta_{ij}=4ij\Delta^{2j}$.
Let $a_i=\frac{1}{2d^2}$. Since $a_i\leq \frac{1}{2}$, we have that
$(1-a_i)\geq \exp(-2a_i)$.
In order to be able to apply the Local Lemma, we need
that for any $1\leq i \leq r$
$$p_i\leq a_i\prod_{j=1}^r \exp(-2a_j\Delta_{ij})\,.$$
That is
$$C^{-i}\leq a^{-i}\prod \exp(-8ija^{-j}d^{2j})\,,$$
or equivalently
$$C\geq a \exp\left( 8\sum_{j=1}^r \frac{j}{2^j}\right)\,.$$
Since the infinite series $\sum_{j=1}^\infty \frac{j}{2^j}$ converges to $2$,
we obtain that for large enough $C$, the conditions of the Local Lemma
are satisfied independently on the size of our finite graph $G$. This ends
the proof of Theorem \ref{alon}.\quad \qed
\vi
Let $|K|=C(|Q|)$ and let $c:\Gamma/H\to K$ be a nonrepetitive $K$-coloring
that gives rise to an element $y\in \sk$. The following
proposition finishes the proof of Theorem \ref{tetel1}.
\begin{propo} \label{propo2}
All elements of the orbit closure $Y$ of $y$ in $\sk$ are $Z$-proper.
\end{propo}
\proof Let $x\in Y$ with underlying Schreier graph $H'$ and
coloring $c':\Gamma/H \to K$. Since $Z$ is an URS, $H'\in Z$. Indeed,
$\pi^{-1}(Z)$ is a closed $\Gamma$-invariant set and $y\in \pi^{-1}(Z)$.
Clearly, $\alpha(\gamma)(x)=x$ if $\gamma\in H'$. Now suppose that
$\alpha(\gamma)(x)=x$ and $\gamma\notin H'$ (that is $x$ is not $Z$-proper).
By Lemma \ref{trivi}, there exists a colored-labeled automorphism $\theta$
of the graph $x$ moving $\roo(x)$ to $\gamma(\roo(x))\neq \roo(x)$. 
Now we proceed similarly as in the proof of Lemma 2. \cite{AGHR} or 
in the proof of Theorem 2.6 \cite{ABT}.
Let $a\in V(x)$ be a vertex such that
there is no $b\in X$ such that
$\dist_x(b,\theta(b))<\dist_x(a,\theta(a))$.
Let $(a=a_1, a_2,\dots, a_{n+1}=\theta(a))$ be a shortest path between $a$ and
$\theta(a)$.
For $1\leq i \leq n$, let
$\gamma_{k_i}(a_i)=a_{i+1}\,.$ Then let $a_{n+2}=\gamma_{k_1}(a_{n+1}),
a_{n+3}=\gamma_{k_1}(a_{n+2}),\dots, a_{2n}=\gamma_{k_n}(a_{2n-1})\,.$
Since $\theta$ is a colored-labeled automorphism, for any $1\leq i \leq n$\,
\begin{equation}\label{e1}
c(a_i)=c(a_{i+n})\,.
\end{equation}
\begin{lemma} 
The walk $(a_1,a_2,\dots, a_{2n})$ is a path.
\end{lemma}
\proof Suppose that the walk above crosses itself, that is for some $i,j$,
$a_j=a_{n+i}$. If $(n+1)-j\geq (n+i)-(n+1)=i-1\,,$ then
$\dist(a_2,\theta(a_2))=\dist(a_2,a_{n+2})<\dist(a,\theta(a))\,.$
On the other hand, if $(n+1)-j\leq (n+i)-(n+1)=i-1\,,$ then
$\dist(a_n,\theta(a_n))=\dist(a_n,a_{2n-1})<\dist(a,\theta(a))\,.$
Therefore, $(a_1,a_2,\dots, a_{2n})$ is a path. \qed
\vi
By (\ref{e1}) and the previous lemma, the $K$-colored Schreier-graph $x$
contains a repetitive path. Since $x$ is in the orbit closure of $y$, this
implies that $y$ contains a repetitive path as well, in contradiction with
our assumption. \qed

\section {Sofic groups and invariant measures}
First, let us recall the notion of a finitely generated sofic group. 
 Let $\G$ be a finitely generated
infinite group with a minimal, symmetric generating system 
$Q=\{\gamma_i\}^r_{i=1}$ and
a surjective homomorphism $\kappa:\F_n\to \Gamma$ from the
free group $\F_n$ with generating system $\overline{Q}=\{r_i\}^n_{i=1}$ mapping
$r_i$ to $\gamma_i$. Let $\Cay^Q_\G$ be the Cayley graph of $\G$ with
respect to the generating system $Q$, that is the Schreier graph corresponding
to the subgroup $H=\{1_\G\}$. Let $\{G_k\}^\infty_{k=1}$ be a sequence
of finite $\F_n$-Schreier graphs.
We call a vertex $p\in V(G_k)$ a $(\Gamma,r)$-vertex if there exists
a rooted isomorphism
$$\Psi:B_r(G_k,p)\to B_r(\Cay^Q_\G,1_\G)\,$$
\vi
such that if $e$ is a directed edge in the ball $B_r(G_k,p)$ labeled
by $r_i$, then the edge $\Phi(e)$ is labeled by $\gamma_i$.
We say that $\{G_k\}^\infty_{k=1}$ is a sofic approximation of $\Cay^Q_\G$, if
for any $r\geq 1$ and a real number $\e>0$ there exists $N_{r,\e}\geq 1$
such that if $k\geq N_{r,\e}$ then there exists a subset $V_k\subset V(G_k)$
consisting of $(\Gamma,r)$-vertices such that $|V_k|\geq (1-\e)|V(G_k)|$.
A finitely generated group $\Gamma$ is called sofic if the Cayley-graphs of $\G$ admit
sofic approximations. Sofic groups were introduced by Gromov in \cite{Gro} under the name
 of initially subamenable groups, the word ``sofic'' was coined by Weiss in
\cite{Weiss}. It is important to note that all the amenable, residually-finite
and residually amenable groups are sofic, but there exist finitely
generated sofic groups that are not residually amenable (see the book
of Capraro and Lupini \cite{Lup} on sofic groups). It is still an open
question whether all groups are sofic.
Now let $\Gamma$ be a finitely generated sofic group with generating system
$Q=\{\gamma_i\}^n_{i=1}$ and a sofic approximation
$\{G_k\}^\infty_{k=1}$. Using Theorem \ref{alon}, for each $k\geq 1$ let us choose a nonrepetitive coloring
$c_k:V(G_k)\to K$, where $|K|\geq C(|Q|)$. We can associate
a probability measure $\mu_k$ on the space
of $K$-colored $\F_n$-Schreier graphs $S^{\overline{Q},K}_{\F_n}$. Note that
the origin of this construction can be traced back to the paper of Benjamini
and Schramm \cite{BS}. For a vertex $p\in V(G_k)$ we consider the
rooted $K$-colored Schreier graph $(G_k^{c_k},p)$. The measure $\mu_k$
is defined as
$$\mu_k=\frac{1}{|V(G_k)|} \sum_{p\in V(G_k)} \delta(G_k^{c_k},p)\,,$$
where $\delta(G_k^{c_k},p)$ is the Dirac-measure on $S^{\overline{Q},K}_{\F_n}$
concentrated on the rooted $K$-colored Schreier graph $(G_k^{c_k},p)$.
Clearly, $\mu_k$ is invariant under the action of $\F_n$.
Since the space of $\F_n$-invariant probability measures on the
compact space $S^{\overline{Q},K}_{\F_n}$ is compact with respect to
the weak-topology, we have a convergent subsequence $\{\mu_{n_k}\}^\infty_{k=1}$
converging weakly to some probability measure $\mu$.
Let $C^{\overline{Q}}_{\F_n}(N)$ be the Schreier graph corresponding to the normal
subgroup $N=\Ker(\kappa)$. This means that we have a natural graph isomorphism
from $C^\Gamma_{\F_n}$ to $\Cay^Q_\G$ that changes the labels $r_i$ to $\gamma_i$.
\begin{propo}
The probability measure $\mu$ is concentrated on the $\F_n$-invariant
closed set $\Omega$ of nonrepetitive $K$-colorings on $C^{\overline{Q}}_{\F_n}(N)$.
\end{propo}
\proof
Let $U_r\subset S^{\overline{Q},K}_{\F_n}$ be the clopen set of
$K$-colored Schreier graphs $G$ such that
the ball $B_r(G,\roo(G))$ is not rooted-labeled isomorphic to 
$B_r(C^{\overline{Q}}_{\F_n}(N),1_\Gamma)$.
By our assumptions on the sofic approximations, $\lim_{k\to\infty}\mu_k(U_r)=0\,,$
hence $\mu(U_r)=0\,.$
Now let $V_r\subset S^{\overline{Q},K}_{\F_n}$ be the clopen set of
$K$-colored Schreier graphs $G$ such that the ball
$B_r(G,\roo(G))$ contains a repetitive path. By our assumptions on the
colorings $c_k$, $\mu_k(V_r)=0$ for any $k\geq 1$. Hence $\mu(V_r)=0$.
Therefore $\mu$ is concentrated on $\Omega$. \qed

\vi
Now we prove Theorem \ref{tetel2}.
Observe that we have an $F_r$-equivariant continuous map $\Sigma:\Omega\to K^\Gamma$, where $\F_n$ acts on the Bernoulli space $K^\Gamma$ on the right by
$\rho(f)(\gamma)=f(\gamma\kappa(\rho))\,$ for $\rho\in\F_n$,$\gamma\in \Gamma$.
Then the image of $F$ is a closed $\Gamma$-invariant subset in $K^\Gamma$, that
is a Bernoulli subshift consisting of elements that are given by nonrepetitive
$K$-colorings. The pushforward of $\mu$ under $\Sigma$ is a $\Gamma$-invariant
probability measure concentrated on $Y$. By Proposition \ref{propo2}, $\Gamma$
acts freely on $Y$, hence Theorem \ref{tetel2} follows. \qed

\vi
\vi
G\'abor Elek

\noindent
Lancaster University

\noindent
g.elek@lancs.ac.uk
\end{document}